\documentclass{article}

\usepackage{PRIMEarxiv}

\usepackage[utf8]{inputenc} 
\usepackage[T1]{fontenc}    
\usepackage{hyperref}       
\usepackage{url}            
\usepackage{booktabs}       
\usepackage{amsfonts}       
\usepackage{nicefrac}       
\usepackage{microtype}      
\usepackage{lipsum}
\usepackage{fancyhdr}       
\usepackage{graphicx}       
\graphicspath{{media/}}     
\usepackage{amsmath}
\usepackage{amsthm}
\usepackage{fdsymbol}
\usepackage{tikz}
\usepackage{float}

\theoremstyle{plain} 
\newtheorem{theorem}{Theorem}[section]
\theoremstyle{definition}
\newtheorem{corollary}[theorem]{Corollary}
\newtheorem{lemma}[theorem]{Lemma}
\newtheorem{proposition}{Proposition}[section]
\newtheorem{definition}{Definition}[section]

\title{Construction of non-regular $A_\alpha$-cospectral graphs from some join of graphs
}

\author{
  Najiya V K, Chithra A V \\
 Department of Mathematics \\
 National Institute of Technology, Calicut\\
 Kerala, India-673601\\
  \texttt{najiya\_p190046ma@nitc.ac.in, chithra@nitc.ac.in} \\
}

\begin{document}
\maketitle

\begin{abstract}
Cospectral graphs are a fascinating concept in graph theory, where two non-isomorphic graphs possess identical sets of eigenvalues. In this paper, we compute the $A_\alpha$-characteristic polynomial of neighbour and non-neighbour splitting join, neighbour and non-neighbour shadow join,  central vertex and edge join and duplicate join of two graphs. In addition, when $G_1$ and $G_2$ are regular, we compute the $A_\alpha$-spectrum of these graphs. As an application, we construct non-regular, non-isomorphic graphs that are $A_\alpha$-cospectral.
\end{abstract}


\section{Introduction}
The graphs in this paper are undirected and simple. There are several matrices associated with a graph, such as adjacency matrix, Laplacian matrix, signless Laplacian matrix etc. Let $G = (V (G) , E (G))$ be a graph with vertex set $V (G) = \{v_1 , v_2 , \dots, v_n \}$ and edge set $E (G)=\{e_1 , e_2 , \dots, e_m \}$. The adjacency matrix $A(G)$ of $G$ is defined by;
$$
(A(G))_{i j}= \begin{cases}1, & \text { if } v_i \text { and } v_j \text { are adjacent; } \\ 0, & \text { otherwise }\end{cases}
$$

The incidence matrix $R(G)$ of a graph $G$ is the matrix of order $n \times m$, whose
$(i, j)$th entry is $1$ if $u_i$ is incident to $e_j$ and $0$ otherwise. $R(G)R(G)^T =
A(G) + D(G)$ and $R(G)R(G)^T = A(G) + rI_n$ for an $r$-regular graph $G$,
where $I_n$ is the identity matrix of order $n$. $R(G)^T R(G) = B(G) + 2I_m $, where $B$ is the adjacency matrix of the line graph of $G$. Let $d_i=d_G\left(v_i\right)$ be the degree of the vertex $v_i$ in $G$, and let $D(G)$ be the diagonal matrix with diagonal entries $d_1, d_2, \ldots, d_n$. The $A_\alpha$-matrix\cite{nikiforov2017merging}, $A_\alpha(G)$, of $G$ is defined as
$A_\alpha(G) =\alpha D (G) +(1-\alpha) A (G)$. It is a convex combination of the degree matrix and adjacency matrix of a graph. For $\alpha=0,\frac{1}{2}$ and $1$, the $A_\alpha$-matrix coincides with the adjacency, signless Laplacian, and degree matrix of $G$. 

For an $n \times n$ matrix $M $, we denote the characteristic polynomial $\det (\lambda I_n - M )$ of $M$ by $\phi(M, \lambda) $. The roots of the $M$-characteristic polynomial of $G$ are the $M$-eigenvalues
of $G$ and the collection of the $M$-eigenvalues, including multiplicities, is called the $M$-spectrum of $G$. If G is an $r$-regular graph, then $\lambda_i(A_\alpha(G))=\alpha r+(1-\alpha)\lambda_i(A(G)).$ Two graphs $G_1$ and $G_2$ are $M$-cospectral if they have the same $M$-spectrum. $A_\alpha$- cospectral graphs are always $A$-cospectral graphs, but $A$-cospectral graphs need not always be $A_\alpha$-cospectral. However, regular $A$-cospectral graphs are $A_\alpha$-cospectral.

The adjacency energy $\varepsilon(G)$\cite{bapat2010graphs} of $G$ is defined as the sum of absolute values of eigenvalues of $G$. The $A_\alpha$-energy $\varepsilon_\alpha(G)$\cite{pirzada2021alpha} of $G$ is defined as $\displaystyle\sum_{i=1}^n\left|\lambda_i(A_\alpha(G))-\frac{2\alpha m}{n}\right|$. For a regular graph $G$, $\varepsilon_\alpha(G)=(1-\alpha)\varepsilon(G)$. In \cite{najiya2023study} introduced the concept of $A_\alpha$-borderenergetic and $A_\alpha$-hyperenergetic graphs. A graph $G$ on $n$ vertices is $A_\alpha$-borderenergetic if $\varepsilon_\alpha(G) = \varepsilon_\alpha(K_n )$, for some $\alpha\in[0,1)$. Borderenergetic graphs are not $A_\alpha$-borderenergetic, but regular borderenergetic graphs are $A_\alpha$-borderenergetic for every value of $\alpha$. The graphs whose $A_\alpha$-energy exceeds the $A_\alpha$-energy of the complete graph on the same vertices are called $A_\alpha$-hyperenergetic. That is, a graph $G$ is $A_\alpha$-hyperenergetic if $\varepsilon_\alpha(G) \geq \varepsilon_\alpha(K_n )$, for some $\alpha\in[0,1)$.

In recent years, numerous researchers have explored the $A_\alpha$-spectral characteristics of graphs created through various graph operations like the Cartesian product, Kronecker product, strong product, lexicographic product, corona, edge corona, and neighbourhood corona, among others. Readers are encouraged to refer to \cite{li2019a_,tahir2020coronae,najiya2023alpha,da2023characteristic,basunia2021a_} and the associated references for a more in-depth examination of these graph operations and the corresponding spectral results. 

Cospectral graphs are a fascinating concept in graph theory, where two non-isomorphic graphs possess identical sets of eigenvalues. Although they have distinct structures, these graphs have the same spectral properties. Cospectral graphs have applications in various fields, including chemistry, computer science, and network analysis. By understanding their properties, we can see how different graphs are alike or different, even if they seem unrelated at first. Construction of non-isomorphic $A_\alpha$-cospectral graphs is a nontrivial problem in spectral graph theory especially for large graphs. There are several methods for finding $A_\alpha$-cospectral families. In this paper, we find infinitely many pairs of non-isomorphic non-regular $A_\alpha$-cospectral graphs using certain join of graphs.

\begin{figure}[H]
\begin{center}
\begin{tikzpicture}  
  [scale=1.1,auto=center] 
 \tikzset{dark/.style={circle,fill=black}}
 \tikzset{hollow/.style={circle,draw=black}}
  \tikzset{white/.style={circle,draw=white}}
    
  \node [dark](a1) at (8,1){} ;  
  \node [dark](a2) at (10,1)  {}; 
  \node [dark] (a3) at (12,1)  {};  
  \node[dark] [dark](a4) at (9,2) {};  
  \node [dark](a5) at (11,2)  {};  
  \node [dark](a6) at (9,3)  {};  
  \node[dark] (a7) at (11,3)  {};  
  \node [dark](a8) at (8,4){};
  \node[dark] (a9) at (10,4) {};
  \node [dark](a10) at (12,4) {};
  
   \node [dark](b1) at (1,1){} ;  
  \node[dark] (b2) at (3,1)  {}; 
  \node [dark](b3) at (5,1)  {};  
  \node [dark](b4) at (2,2) {};  
  \node[dark] (b5) at (4,2)  {};  
  \node [dark](b6) at (2,3)  {};  
  \node [dark](b7) at (4,3)  {};  
  \node[dark] (b8) at (1,4){};
  \node [dark](b9) at (3,4) {};
  \node [dark](b10) at (5,4) {};
  
    \node  at (3,0) {(a) $G_1$};
    \node  at (10,0) {(b) $G_2$};

  \draw (b1) -- (b2);
  \draw (b2) -- (b3);  
  \draw (b2) -- (b4);  
  \draw (b4) -- (b6);  
  \draw (b3) -- (b4);  
  \draw (b6) -- (b8);  
  \draw (b5) -- (b7);
  \draw (b2) -- (b5);  
  \draw (b8) -- (b4);  
  \draw (b9) -- (b6);  
  \draw (b3) -- (b7);  
  \draw (b9) -- (b7);  
  \draw (b5) -- (b1);
  \draw (b10) -- (b3);  
  \draw (b1) -- (b8);  
  \draw (b1) -- (b6);  
  \draw (b10) -- (b5);  
  \draw (b10) -- (b7);  
  \draw (b8) -- (b9);
  \draw (b10) -- (b9);
  
  \draw (a1) -- (a2);
  \draw (a2) -- (a3);  
  \draw (a2) -- (a4);  
  \draw (a4) -- (a6);  
  \draw (a3) -- (a7);  
  \draw (a6) -- (a8);  
  \draw (a5) -- (a7);
  \draw (a2) -- (a5);  
  \draw (a8) -- (a5);  
  \draw (a9) -- (a6);  
  \draw (a3) -- (a5);  
  \draw (a9) -- (a7);  
  \draw (a4) -- (a1); 
  \draw (a10) -- (a3);  
  \draw (a1) -- (a8);  
  \draw (a1) -- (a6);  
  \draw (a10) -- (a4);  
  \draw (a10) -- (a7);  
  \draw (a8) -- (a9);
  \draw (a10) -- (a9);
  
\end{tikzpicture}  

\end{center}
\caption{Non-isomorphic $A$-cospectral regular graphs} \label{acoreg}
\end{figure}
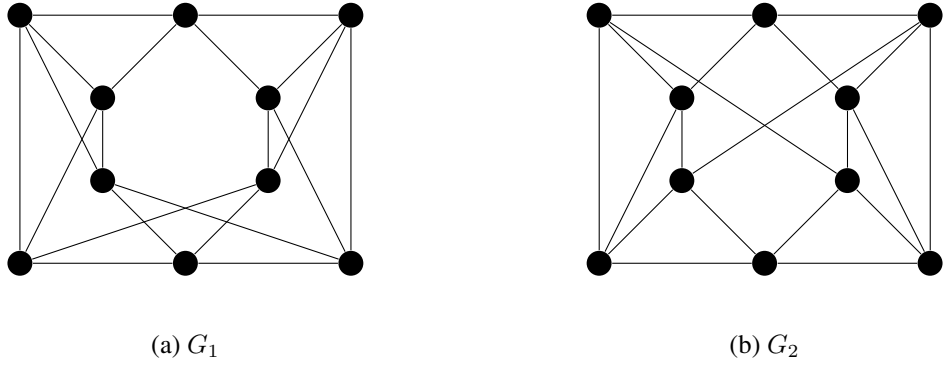

We will use the symbols $O_{m\times n}$ and $J_{m\times n}$ for the $m \times n$
matrices consisting of all $0$’s and all $1$’s, respectively. For convenience, the adjacency matrix of $G_i$, $A(G_i)$ will be denoted as $A_i$, its degree matrix as $D_i$ and the $A_\alpha$-matrix as $A_{\alpha_i}$. $\overline{A_i}$ denotes the adjacency matrix of complement of $G$, which is a graph $\overline{G}$ on the same set of vertices as $G$ such that there will be an edge between two vertices $(v_i, v_j)$ in $\overline{G}$, if and only if there is no edge between $(v_i, v_j)$ in $G$.

The structure of the paper is as follows. After preliminaries in Section \ref{prelim}, we compute the $A_\alpha $-characteristic polynomial and $A_\alpha$-spectrum of different joins of graphs and use the results to construct cospectral graphs in Section \ref{main}.

\section{Preliminaries}\label{prelim}
\begin{definition}\cite{balakrishnan2012textbook}
    Let $G_1$ and $G_2$ be two vertex disjoint graphs. Then the join $G_1 \vee G_2$ of $G_1$ and $G_2$ is the graph in which each vertex of $G_1$ is made adjacent to every vertex of $G_2$.
\end{definition}

\begin{definition}\cite{lu2023spectra}
    Let $G_1$ and $G_2$ be two vertex disjoint graphs with $V\left(G_1\right)=\left\{u_1, u_2, \ldots, u_n\right\}$. The neighbour splitting(V-vertex) join of $G_1$ and $G_2$, denoted by $G_1 \barwedge G_2$, is obtained by adding vertices $u_1^{\prime}, u_2^{\prime}, \ldots, u_n^{\prime}$ to $G_1 \vee G_2$ and connecting $u_i^{\prime}$ to $u_j$ if and only if $\left(u_i, u_j\right) \in E\left(G_1\right)$.
\end{definition}

\begin{definition}\cite{hamud2023new}
    Let $G_1$ and $G_2$ be two vertex disjoint graphs with $V\left(G_1\right)=\left\{u_1, u_2, \ldots, u_n\right\}$. The non-neighbour splitting join of $G_1$ and $G_2$, denoted by $G_1 \doublebarwedge G_2$, is obtained by adding the vertices $u_1^{\prime}, u_2^{\prime}, \ldots, u_n^{\prime}$ to $G_1 \vee G_2$ and connecting $u_i^{\prime}$ to $u_j$ if and only if $\left(u_i, u_j\right) \notin E\left(G_1\right)$.
\end{definition}

\begin{definition}\cite{jahfar2020central}
    The central vertex join of $G_1$ and $G_2$, denoted by $G_1\veedot G_2$, is the graph obtained by subdividing each edge of $G_1$ exactly once and joining all non-adjacent vertices in $G_1$ and then joining each vertex of $G_1$ to every vertex of $G_2$.
\end{definition}

\begin{definition}\cite{jahfar2020central}
    The central edge join of $G_1$ and $G_2$, denoted by $G_1\veebar G_2$, is the graph obtained by subdividing each edge of $G_1$ exactly once and joining all the non-adjacent vertices in $G_1$ and then joining each vertex corresponding to edges of $G_1$ to every vertex of $G_2$.
\end{definition}

\begin{definition}\cite{sampathkumar1973duplicate}
    Let $G=(V, E)$ be a simple graph. Let $V^{\prime}$ be a set such that $|V|=\left|V^{\prime}\right|, V \cap V^{\prime}=\emptyset$ and $f: V \rightarrow V^{\prime}$ be bijective(for $v\in V$ we write $f(v)=v'$). A duplicate graph of $G$ is $D(G)=\left(V_1, E_1\right)$, where the set of vertices $V_1=V \cup V^{\prime}$ and the set of edges $E_1$ of $D(G)$ is defined as, the edges $uv^{\prime}$ and $u^{\prime} v$ are in $E_1$ if and only if $uv$ is in $E$.
\end{definition}

Inspired by these operations, we define new joins of graphs such as neighbour and non-neighbour shadow join and duplicate join of graphs.

\begin{definition}
    Let $G_1$ and $G_2$ be two vertex disjoint graphs with $V\left(G_1\right)=\left\{u_1, u_2, \ldots, u_n\right\}$. The neighbour shadow join of $G_1$ and $G_2$, denoted by $G_1 \doublesqcup G_2$, is obtained by taking a copy of $G_1$, say $G_1'$, to $G_1 \vee G_2$ and connecting $u_i^{\prime}$ to $u_j$ if and only if $\left(u_i, u_j\right) \in E\left(G_1\right)$.
\end{definition}

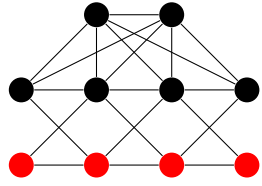
\begin{figure}[H]
    \centering
    \begin{tikzpicture}[scale=1,auto=center] 
 \tikzset{dark/.style={circle,fill=black}}
 \tikzset{red/.style={circle,fill=red}}
 \tikzset{white/.style={circle,draw=white}}
  
 \node [dark] (a1) at (1,1)  {};  
  \node [dark] (a2) at (2,1)  {};
  
  \node [dark] (a3) at (0,0)  {};  
  \node [dark] (a4) at (1,0){};
  \node [dark] (a5) at (2,0) {};
  \node [dark] (a6) at (3,0) {};
  
   \node [red] (a7) at (0,-1)  {};  
  \node [red] (a8) at (1,-1){};
  \node [red] (a9) at (2,-1) {};
  \node [red] (a10) at (3,-1) {};

  \draw (a1) -- (a2);  
  
  \draw (a3) -- (a4);  
  \draw (a4) -- (a5);
  \draw (a5) -- (a6);  
  
  \draw (a1) -- (a3);  
  \draw (a1) -- (a4);  
  \draw (a1) -- (a5);  
  \draw (a1) -- (a6);  

   \draw (a2) -- (a3);  
  \draw (a2) -- (a4);  
  \draw (a2) -- (a5);  
  \draw (a2) -- (a6);
  
  \draw (a7) -- (a8);   
  \draw (a8) -- (a9);
  \draw (a9) -- (a10);

  \draw (a7) -- (a4);   
  \draw (a8) -- (a3);
  \draw (a8) -- (a5);   
  \draw (a9) -- (a4);
  \draw (a9) -- (a6);   
  \draw (a10) -- (a5);

\end{tikzpicture}
    \caption{$P_4\doublesqcup P_2$}
    \label{neighshadow}
\end{figure}

\begin{definition}
    Let $G_1$ and $G_2$ be two vertex disjoint graphs with $V\left(G_1\right)=\left\{u_1, u_2, \ldots, u_n\right\}$. The non neighbours shadow join of $G_1$ and $G_2$, denoted by $G_1 \doublesqcap G_2$, is obtained by taking a copy of $G_1$, say $G_1'$, to $G_1 \vee G_2$ and connecting $u_i^{\prime}$ to $u_j$ if and only if $\left(u_i, u_j\right) \notin E\left(G_1\right)$.
\end{definition}

\begin{figure}[H]
    \centering
    \begin{tikzpicture}[scale=1,auto=center] 
 \tikzset{dark/.style={circle,fill=black}}
 \tikzset{red/.style={circle,fill=red}}
 \tikzset{white/.style={circle,draw=white}}
  
 \node [dark] (a1) at (1,1)  {};  
  \node [dark] (a2) at (2,1)  {};
  
  \node [dark] (a3) at (0,0)  {};  
  \node [dark] (a4) at (1,0){};
  \node [dark] (a5) at (2,0) {};
  \node [dark] (a6) at (3,0) {};
  
   \node [red] (a7) at (0,-1)  {};  
  \node [red] (a8) at (1,-1){};
  \node [red] (a9) at (2,-1) {};
  \node [red] (a10) at (3,-1) {};

  \draw (a1) -- (a2);  
  
  \draw (a3) -- (a4);  
  \draw (a4) -- (a5);
  \draw (a5) -- (a6);  
  
  \draw (a1) -- (a3);  
  \draw (a1) -- (a4);  
  \draw (a1) -- (a5);  
  \draw (a1) -- (a6);  

   \draw (a2) -- (a3);  
  \draw (a2) -- (a4);  
  \draw (a2) -- (a5);  
  \draw (a2) -- (a6);
  
  \draw (a7) -- (a8);   
  \draw (a8) -- (a9);
  \draw (a9) -- (a10);

  \draw (a7) -- (a5);   
  \draw (a7) -- (a6);
  \draw (a8) -- (a6);   
  \draw (a9) -- (a3);
  \draw (a10) -- (a3);   
  \draw (a10) -- (a4);

\end{tikzpicture}
    \caption{$P_4\doublesqcap P_2$}
    \label{nonneighshadow}
\end{figure}
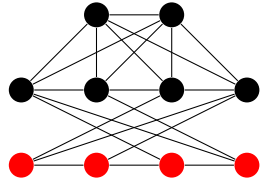

\begin{definition}
    Let $G_1$ and $G_2$ be any two graphs on $n_1 , n_2$ vertices and $m_1 , m_2$ edges, respectively. The duplicate join of $G_1$ and $G_2$ is the graph $G_1\bowtie G_2$, obtained from $D(G_1 )$ and $G_2$ by joining each vertex in $V(G_1)$ with every vertex of $G_2$.
\end{definition}

\begin{figure}[H]
    \centering
    \begin{tikzpicture}[scale=1,auto=center] 
 \tikzset{dark/.style={circle,fill=black}}
 \tikzset{red/.style={circle,fill=red}}
 \tikzset{white/.style={circle,draw=white}}
  
 \node [dark] (a1) at (1,1)  {};  
  \node [dark] (a2) at (2,1)  {};
  
  \node [dark] (a3) at (0,0)  {};  
  \node [dark] (a4) at (1,0){};
  \node [dark] (a5) at (2,0) {};
  \node [dark] (a6) at (3,0) {};
  
   \node [red] (a7) at (0,-1)  {};  
  \node [red] (a8) at (1,-1){};
  \node [red] (a9) at (2,-1) {};
  \node [red] (a10) at (3,-1) {};

  \draw (a1) -- (a2);  
  
  
  \draw (a1) -- (a3);  
  \draw (a1) -- (a4);  
  \draw (a1) -- (a5);  
  \draw (a1) -- (a6);  

   \draw (a2) -- (a3);  
  \draw (a2) -- (a4);  
  \draw (a2) -- (a5);  
  \draw (a2) -- (a6);
  

  \draw (a7) -- (a4);   
  \draw (a8) -- (a3);
  \draw (a8) -- (a5);   
  \draw (a9) -- (a4);
  \draw (a9) -- (a6);   
  \draw (a10) -- (a5);

\end{tikzpicture}
    \caption{$P_4\bowtie P_2$}
    \label{duplicate}
\end{figure}
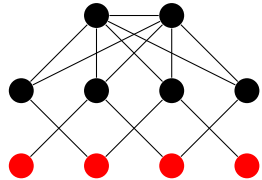

\begin{definition}\cite{zhang2006schur}
    Let $M$ be a block matrix $$M=\begin{bmatrix}
        A & B\\
        C & D
    \end{bmatrix}$$ 

such that its blocks $A$ and $D$ are square. If $A$ is invertible, the Schur complement of $A$ in $M$ is
$$
M / A=D-C A^{-1} B
$$
and if $D$ is invertible, the Schur complement of $D$ in $M$ is
$$
M / D=A-B D^{-1} C .
$$
\end{definition}

\begin{lemma}\cite{zhang2006schur}\label{schur}
    If $D$ is invertible then,
$$
|M|=|M / D||D|
$$
and if $A$ is invertible then,
$$
|M|=|M / A||A|
$$
\end{lemma}

\begin{lemma}\cite{hamud2023new}\label{schurcomp}
    Let $M$ be a block matrix
$$
M=\left(\begin{array}{ccc}
A & B & J_{n_1 \times n_2} \\
B & C & O_{n_1 \times n_2} \\
J_{n_2 \times n_1} & O_{n_2 \times n_1} & D
\end{array}\right)
$$
where $A, B$, and $C$ are square matrices of order $n_1$ and $D$ is a square matrix of order $n_2$. Then the Schur complement of $x I_{n_2}-D$ in the characteristic matrix of $M$ is
$$
\left(\begin{array}{cc}
x I_{n_1}-A-\Gamma_D(x) J_{n_1} & -B \\
-B & x I_{n_1}-C
\end{array}\right)
$$
\end{lemma}

\begin{lemma}\cite{cvetkovic1980spectra}\label{pa}
	Let $G$ be an $r$-regular connected graph of order $n$ with adjacency matrix $A$ having $t$ distinct eigenvalues $r=\lambda_1,\lambda_2,...,\lambda_t$.  Then there exists a polynomial \[P(\lambda)=n\dfrac{(\lambda-\lambda_2)(\lambda-\lambda_3)...(\lambda-\lambda_t)}{(r-\lambda_2)(r-\lambda_3)...(r-\lambda_t)}.\] such that $P(A)=J_n, P(r)=n$ and $P(\lambda_i) =0$ for $\lambda_i\neq r.$
\end{lemma}

\begin{definition}\cite{cui2012spectrum}\label{gammad}
	The $M$-coronal $\Gamma_M(x)$  of a matrix $M$ of order $n$ is defined as the sum of the entries in the matrix $(xI_n-M)^{-1}$ $($if it exists$)$, that is, $$\Gamma_M(x)=J_{n\times 1}^T(xI_n-M)^{-1}J_{n\times 1}.$$ \\
	If each row sum of an $n\times n$ matrix $M$ is constant, say a, then 
	$\Gamma_{M}(x)=\frac{n}{x-a}$.
\end{definition}

\begin{lemma}\cite{das2019new}\label{inverse}
    For any real numbers $c, d > 0$, $(cI_n - dJ_n )^{-1} = \frac{1}{c} I_n + \frac{c}{c(c-nd)}J_n$.
\end{lemma}

\section{Main results}\label{main}
Here we compute the characteristic polynomials of neighbour and non-neighbour splitting join, neighbour and non-neighbour shadow join and duplicate join of two arbitrary graphs. Also, we compute the characteristic polynomial of the central vertex and the edge joining of $G_1 $ and $ G_2$ where $G_1$ is regular. Then we compute the $A_\alpha$-spectrum of these operations when $G_1$ and
$G_2$ are regular. As an application, we then construct infinitely many pairs of non-regular $A_\alpha$-cospectral graphs.

Throughout the section $A_i$, $A_{\alpha_i}$, $\overline{A_i}$ and $D_i$ represents $A(G_i)$, $A_{\alpha}(G_i)$, $A(\overline{G_i})$ and $D(G_i)$ respectively.

Before we dive into our results, we introduce a useful lemma derived for matrices with a certain structure. This lemma simplifies the computation of their characteristic polynomial, a major step in our results. 
\begin{lemma}\label{schurmy}
     Let $M$ be a block matrix
$$
M=\left(\begin{array}{ccc}
A & B & J_{n_1 \times n_2} \\
B & C & O_{n_1 \times n_2} \\
J_{n_2 \times n_1} & O_{n_2 \times n_1} & D
\end{array}\right)
$$
where $A, B$, and $C$ are square matrices of order $n_1$ and $D$ is a square matrix of order $n_2$. Then the characteristic polynomial of $M$ is 
$$\left|\lambda I_{n_2}-D\right|\left|\lambda I_{n_1}-C\right|\left|\lambda I_{n_1}-A-B(\lambda I_{n_1}-C)^{-1}B\right|\left[1-\Gamma_D(\lambda)\Gamma_{A+B(\lambda I_{n_1}-C)^{-1}B}(\lambda)\right]$$
\end{lemma}

\begin{proof}
    \begin{align*}
    \intertext{Given}
    M=&\left(\begin{array}{ccc}
        A & B & J \\
        B & C & 0 \\
        J & 0 & D
    \end{array}\right).\\
    \intertext{Then, the characteristic polynomial of $M$ is}
        \phi(M,\lambda)=&\left|\begin{array}{ccc}
           \lambda-A  & -B & -J \\
            -B & \lambda-C & 0\\
            -J & 0 & \lambda-D
        \end{array}\right|.\\
        \intertext{By Lemmas \ref{schur} and \ref{schurcomp},}
        \phi(M,\lambda)=&\left|\lambda I-D\right|\left|\begin{array}{cc}
            \lambda I-A-\Gamma_D(\lambda)J & -B \\
            -B & \lambda I-C
        \end{array}\right|\\
        =&\left|\lambda I-D\right|\left|\lambda I-C\right|\left|\lambda I-A-\Gamma_D(\lambda)J-B(\lambda I-C)^{-1}B\right|\\
        =&\left|\lambda I-D\right|\left|\lambda I-C\right|\left[\left|\lambda I-A-B(\lambda I-C)^{-1}B\right|\left[\Gamma_D(\lambda)\mathbf{1_n^T}adj(\lambda I-A-B(\lambda I-C)^{-1}B)\mathbf{1_n}\right]\right]\\
        =&\left|\lambda I-D\right|\left|\lambda I-C\right|\left|\lambda I-A-B(\lambda I-C)^{-1}B\right|\left[1-\Gamma_D(\lambda)\Gamma_{A+B(\lambda I-C)^{-1}B}(\lambda)\right].
    \end{align*}
\end{proof}
\subsection{Splitting joins}
In this section, we derive the $A_\alpha$-characteristic polynomial of $G_1 \barwedge G_2$ and $G_1 \doublebarwedge G_2$ when $G_1$ and $G_2$ are arbitrary graphs.  

\begin{proposition}\label{split}
    Let $G_i$ be a graph on $n_i$ vertices for $i=1,2$. Then

\begin{align*}
\phi(A_\alpha  \left(G_1 \barwedge G_2\right),\lambda) =&\left|(\lambda-\alpha n_1)I-A_{\alpha_2}\right| \left|\lambda I-\alpha D_1\right| \\
& \left|(\lambda-\alpha n_2)I-\alpha D_1-A_{\alpha_1}-(1-\alpha)^2 A_1(\lambda I-\alpha D_1)^{-1}A_1\right| \\
& \left[1-\Gamma_{A_{\alpha_2}}\left(\lambda-\alpha n_1\right) \Gamma_{\alpha D_1+A_{\alpha_1}+(1-\alpha)^2A_1(\lambda I-\alpha D_1)^{-1}A_1}\left(\lambda-\alpha n_2\right)\right] .\\
\\
\phi(A_\alpha  \left(G_1 \doublebarwedge G_2\right),\lambda) =&\left|(\lambda-\alpha n_1)I-A_{\alpha_2}\right| \left|(\lambda-\alpha(n-1))I+\alpha D_1\right| \\
&  \left|(\lambda-\alpha(n_1+n_2-1))I-(1-\alpha)A_1-(1-\alpha)^2\overline{A_1}((\lambda-\alpha(n_1-1))I+\alpha D_1)^{-1}\overline{A_1}\right| \\
& \left[1-\Gamma_{A_{\alpha_2}}\left(\lambda-\alpha n_1\right) \Gamma_{(1-\alpha)A_1+(1-\alpha)^2\overline{A_1}((\lambda-\alpha(n_1-1))I+\alpha D_1)^{-1}\overline{A_1}}\left(\lambda-\alpha(n_1+n_2-1)\right)\right] .
\end{align*}

\end{proposition}

\begin{proof}
   With a suitable ordering of the vertices of $G_1 \barwedge G_2$ and $G_1 \doublebarwedge G_2$, we get

\begin{align*}
A_\alpha\left(G_1 \barwedge G_2\right)=&\left(\begin{array}{ccc}
\alpha n_2 I_{n_1}+\alpha D_1+A_{\alpha_1} & (1-\alpha)A_1 & (1-\alpha)J_{n_1 \times n_2} \\
(1-\alpha)A_1 & \alpha D_1 & O_{n_1 \times n_2} \\
(1-\alpha)J_{n_2 \times n_1} & O_{n_2 \times n_1} & \alpha n_1 I_{n_2}+A_{\alpha_2}
\end{array}\right),\\
\intertext{and}
 A_\alpha\left(G_1 \doublebarwedge G_2\right)=&\left(\begin{array}{ccc}
\alpha\left(n_1+n_2-1\right) I_{n_1}+(1-\alpha)A_1 & (1-\alpha)\overline{A_1} & (1-\alpha)J_{n_1 \times n_2} \\
(1-\alpha)\overline{A_1} & \alpha\left(\left(n_1-1\right) I_{n_1}-D_1\right) & O_{n_1 \times n_2} \\
(1-\alpha)J_{n_2 \times n_1} & O_{n_2 \times n_1} & \alpha\left(D_2+n_1 I_{n_2}\right)+(1-\alpha)A_2
\end{array}\right).
\end{align*}
Then the results follow from Lemma \ref{schurmy}.

\end{proof}

Now, in the following propositions, we obtain the $A_\alpha $-eigenvalues of $G_1 \barwedge G_2$ and $G_1 \doublebarwedge G_2$ where $G_i$'s are $r_i$-regular graphs.
\begin{proposition}\label{spectrum_nspj}
    Let $G_i$ be $r_i$-regular graph with $n_i$ vertices for $i=1,2$. Then the $A_\alpha$-spectrum of $G_1 \barwedge G_2$ consists of:
    \begin{enumerate}
        \item $\alpha(n_1-r_2)+(1-\alpha)\lambda_i(A_2)$ for each $i=2,3,\dots,n_2$,
        \item two roots of the equation $$(\lambda-\alpha r_1)(\lambda-\alpha(n_2+2r_1)-(1-\alpha)\lambda_i(A_1))-(1-\alpha)^2\lambda_i(A_1)^2=0, \text{ for each $i=2,3\dots,n_1$,}$$
        \item three roots of the equation
        $$(\lambda-\alpha n_1-r_2)((\lambda-\alpha r_1)(\lambda-\alpha n_2-(1+\alpha)r_1)-(1-\alpha)^2r_1^2)-n_1n_2(\lambda-\alpha r_1)=0.$$
    \end{enumerate}
\end{proposition}

\begin{proposition}\label{spectrum_nnspj}
    Let $G_i$ be $r_i$-regular graph with $n_i$ vertices for $i=1,2$. Then the $A_\alpha$-spectrum of $G_1 \doublebarwedge G_2$ consists of:
    \begin{enumerate}
        \item $\alpha(n_1+r_1)-(1-\alpha)\lambda_i(A_2)$ for each $i=2,3,\dots,n_2$,
        \item two roots of the equation $$(\lambda-\alpha(n_1-r_1-1))(\lambda-\alpha(n_1+n_2-1)-(1-\alpha)\lambda_i(A_1))+(1-\alpha)^2(1+\lambda_i(A_1))=0, \text{ for each $i=2,3\dots,n_1$,}$$
        \item three roots of the equation
        $$(\lambda-\alpha n_1-r_2)((\lambda-\alpha(n_1-r_1-1))(\lambda-\alpha(n_1+n_2-1)-(1-\alpha)r_1)-(1-\alpha)^2(n_1-1-r_1)^2)-n_1n_2(\lambda-\alpha(n_1-1)+\alpha r_1)=0.$$
    \end{enumerate}
\end{proposition}

Using Propositions \ref{spectrum_nspj} and \ref{spectrum_nnspj} we generate new families of non-isomorphic non-regular $A_\alpha $-cospectral graphs in the following corollaries. These can be used to obtain $A$-cospectral and $Q$-cospectral graphs.
\begin{corollary}

    \begin{enumerate}
    \item Let $G_1$ and $G_2$ be two $A_\alpha$-cospectral regular graphs and $H$ be an arbitrary graph. Then the graphs $G_1\barwedge H$ and $G_2\barwedge H$ are $A_\alpha$-cospectral.
    \item Let $H_1$ and $H_2$ be two $A_\alpha$-cospectral graphs with $\Gamma_{A_\alpha(H_1)}(\lambda-\alpha n_1) = \Gamma_{A_\alpha(H_2)}(\lambda-\alpha n_1)$ for $\alpha \in [0, 1]$. If $G$ is a regular graph, then the graphs $G\barwedge H_1$ and $G\barwedge H_2$ are $A_\alpha$-cospectral.
\end{enumerate}
\end{corollary}

\begin{corollary}

    \begin{enumerate}
    \item Let $G_1$ and $G_2$ be two $A$-cospectral regular graphs and $H$ be an arbitrary graph. Then the graphs $G_1\doublebarwedge H$ and $G_2\doublebarwedge H$ are $A_\alpha$-cospectral.
    \item Let $H_1$ and $H_2$ be two $A_\alpha$-cospectral graphs with $\Gamma_{A_\alpha(H_1)}(\lambda-\alpha n_1) = \Gamma_{A_\alpha(H_2)}(\lambda-\alpha n_1)$ for $\alpha \in [0, 1]$. If $G$ is a regular graph, then the graphs $G\doublebarwedge H_1$ and $G\doublebarwedge H_2$ are $A_\alpha$-cospectral.
\end{enumerate}
\end{corollary}

\subsection{Shadow joins}
In this section, we derive the $A_\alpha$-characteristic polynomial of $G_1 \doublesqcup G_2$ and $G_1 \doublesqcap G_2$ when $G_1$ and $G_2$ are arbitrary graphs.
\begin{proposition}\label{shadow}
     Let $G_i$ be a graph on $n_i$ vertices for $i=1,2$. Then

\begin{align*}
\phi(A_\alpha  \left(G_1 \doublesqcup G_2\right),\lambda) =&\left|(\lambda-\alpha n_1)I-A_{\alpha_2}\right| \left|\lambda I-\alpha D_1-A_{\alpha_1}\right| \\
&  \left|(\lambda-\alpha n_2)I-\alpha D_1-A_{\alpha_1}-(1-\alpha)^2\overline{A_1}(\lambda I-\alpha D_1-A_{\alpha_1})^{-1}\overline{A_1}\right| \\
& \left[1-\Gamma_{A_{\alpha_2}}\left(\lambda-\alpha n_1\right) \Gamma_{\alpha D_1+A_{\alpha_1}+(1-\alpha)^2\overline{A_1}(\lambda I-\alpha D_1-A_{\alpha_1})^{-1}\overline{A_1}}\left(\lambda-\alpha n_2\right)\right] .\\
\\
\phi(A_\alpha  \left(G_1 \doublesqcap G_2\right),\lambda) =&\left|(\lambda-\alpha n_1)I-A_{\alpha_2}\right| \left|(\lambda-\alpha(n_1-1))I+(1-\alpha) A_1\right| \\
&  \left|(\lambda-\alpha(n_1+n_2-1))I-(1-\alpha)A_1-(1-\alpha)^2\overline{A_1}((\lambda-\alpha(n_1-1))I+(1-\alpha) A_1)^{-1}\overline{A_1}\right| \\
& \left[1-\Gamma_{A_{\alpha_2}}\left(\lambda-\alpha n_1\right) \Gamma_{(1-\alpha)A_1+(1-\alpha)^2\overline{A_1}((\lambda-\alpha(n_1-1))I+(1-\alpha) A_1)^{-1}\overline{A_1}}\left(\lambda-\alpha(n_1+n_2-1)\right)\right] .
\end{align*}

\end{proposition}
\begin{proof}
   With a suitable ordering of the vertices of $G_1 \doublesqcup G_2$ and $G_1 \doublesqcap G_2$, we get

\begin{align*}
A_\alpha\left(G_1 \doublesqcup G_2\right)=&\left(\begin{array}{ccc}
\alpha n_2 I_{n_1}+\alpha D_1+A_{\alpha_1} & (1-\alpha)A_1 & (1-\alpha)J_{n_1 \times n_2} \\
(1-\alpha)A_1 & \alpha D_1+A_{\alpha_1} & O_{n_1 \times n_2} \\
(1-\alpha)J_{n_2 \times n_1} & O_{n_2 \times n_1} & \alpha n_1 I_{n_2}+A_{\alpha_2}
\end{array}\right),\\
\intertext{and}
 A_\alpha\left(G_1 \doublesqcap G_2\right)=&\left(\begin{array}{ccc}
\alpha\left(n_1+n_2-1\right) I_{n_1}+(1-\alpha)A_1 & (1-\alpha)\overline{A_1} & (1-\alpha)J_{n_1 \times n_2} \\
(1-\alpha)\overline{A_1} & \alpha\left(n_1-1\right) I_{n_1}+(1-\alpha)A_1 & O_{n_1 \times n_2} \\
(1-\alpha)J_{n_2 \times n_1} & O_{n_2 \times n_1} & \alpha n_1 I_{n_2}+A_{\alpha_2}
\end{array}\right).
\end{align*}
Then the results follow from Lemma \ref{schurmy}.

\end{proof}

Now, in the following propositions, we obtain the $A_\alpha $-eigenvalues of $G_1 \doublesqcup G_2$ and $G_1 \doublesqcap G_2$ where $G_i$'s are $r_i$-regular graphs.

\begin{proposition}
    Let $G_i$ be $r_i$-regular graph with $n_i$ vertices for $i=1,2$. Then the $A_\alpha$-spectrum of $G_1 \doublesqcup G_2$ consists of:
    \begin{enumerate}
        \item $\alpha(n_1+r_2)+(1-\alpha)\lambda_i(A_2)$ for each $i=2,3,\dots,n_2$,
        \item two roots of the equation $$(\lambda-2\alpha r_1-(1-\alpha)\lambda_i(A_1))(\lambda-\alpha(n_2+2r_1)-(1-\alpha)\lambda_i(A_1))-(1-\alpha)^2\lambda_i(A_1)^2=0, \text{ for each $i=2,3\dots,n_1$,}$$
        \item three roots of the equation
        $$(\lambda-\alpha n_1-r_2)((\lambda- r_1(1+\alpha))(\lambda-\alpha n_2-(1+\alpha)r_1)-(1-\alpha)^2r_1^2)-n_1n_2(\lambda-r_1(1+\alpha))=0.$$
    \end{enumerate}
\end{proposition}
\begin{proposition}
    Let $G_i$ be $r_i$-regular graph with $n_i$ vertices for $i=1,2$. Then the $A_\alpha$-spectrum of $G_1 \doublesqcap G_2$ consists of:
    \begin{enumerate}
        \item $\alpha(n_1+r_2)+(1-\alpha)\lambda_i(A_2)$ for each $i=2,3,\dots,n_2$,
        \item two roots of the equation $$(\lambda-\alpha(n_1-1)+(1-\alpha)\lambda_i(A_1))(\lambda-\alpha(n_1+n_2-1)-(1-\alpha)\lambda_i(A_1))+(1-\alpha)^2(1+\lambda_i(A_1))^2=0, \text{ for each $i=2,3\dots,n_1$,}$$
        \item three roots of the equation\\
        $(\lambda-\alpha n_1-r_2)((\lambda-\alpha(n_1-1)+(1-\alpha)r_1)(\lambda-\alpha(n_1+n_2-1)-(1-\alpha)r_1)-(1-\alpha)^2(n_1-1-r_1)^2)-n_1n_2(\lambda-\alpha(n_1-1)+(1-\alpha) r_1)=0.$
    \end{enumerate}
\end{proposition}

In the upcoming corollaries, we generate new pairs of graphs that are $A_\alpha$-cospectral.

\begin{corollary}\label{cosnshj}

    \begin{enumerate}
    \item Let $G_1$ and $G_2$ be two $A_\alpha$-cospectral regular graphs, and $H$ be an arbitrary graph. Then the graphs $G_1\doublesqcup H$ and $G_2\doublesqcup H$ are $A_\alpha$-cospectral.
    \item Let $H_1$ and $H_2$ be two $A_\alpha$-cospectral graphs with $\Gamma_{A_\alpha(H_1)}(\lambda-\alpha n_1) = \Gamma_{A_\alpha(H_2)}(\lambda-\alpha n_1)$ for $\alpha \in [0, 1]$. If $G$ is a regular graph, then the graphs $G\doublesqcup H_1$ and $G\doublesqcup H_2$ are $A_\alpha$-cospectral.
\end{enumerate}
\end{corollary}

\begin{figure}[H]
\begin{center}
\begin{tikzpicture}[scale=.5,auto=center] 
 \tikzset{dark/.style={circle,fill=black}}
 \tikzset{hollow/.style={circle,draw=black}}
 \tikzset{white/.style={circle,draw=white}}
    \node [dark] (a1) at (-3.5,-3){} ;  
  \node [dark] (a2) at (0,-4)  {}; 
  \node [dark] (a3) at (3.5,-3)  {};  
  \node [dark] (a4) at (-1.5,-1) {};  
  \node [dark] (a5) at (1.5,-1)  {};  
  \node [dark] (a6) at (-1.5,1)  {};  
  \node [dark] (a7) at (1.5,1)  {};  
  \node [dark] (a8) at (-3.5,3){};
  \node [dark] (a9) at (0,4) {};
  \node [dark] (a10) at (3.5,3) {};

  \node [hollow] (b1) at (-5.5,0){} ;  
   \node [hollow] (b2) at (5.5,0){} ;  

   \node [hollow] (b3) at (-5.5,-5){} ;  
   \node [hollow] (b4) at (5.5,-5){} ;  
  
  \node  at (0,-7) {(c) $K_2\doublesqcup G_1$};
   \node  at (13,-7) {(d) $K_2\doublesqcup G_2$};
  
  \draw[red, very thick] (a1) -- (a2);
  \draw[red, very thick] (a2) -- (a3);  
  \draw[red, very thick] (a2) -- (a4);  
  \draw[red, very thick] (a4) -- (a6);  
  \draw[red, very thick] (a3) -- (a7);  
  \draw[red, very thick] (a6) -- (a8);  
  \draw[red, very thick] (a5) -- (a7);
  \draw[red, very thick] (a2) -- (a5);  
  \draw[red, very thick] (a8) -- (a5);  
  \draw[red, very thick] (a9) -- (a6);  
  \draw[red, very thick] (a3) -- (a5);  
  \draw[red, very thick] (a9) -- (a7);  
  \draw[red, very thick] (a4) -- (a1); 
  \draw[red, very thick] (a10) -- (a3);  
  \draw[red, very thick] (a1) -- (a8);  
  \draw[red, very thick] (a1) -- (a6);  
  \draw[red, very thick] (a10) -- (a4);  
  \draw[red, very thick] (a10) -- (a7);  
  \draw[red, very thick] (a8) -- (a9);
  \draw[red, very thick] (a10) -- (a9);
  
  \draw (b1) -- (a2);
 \draw (b1) -- (a2);
 \draw (b1) -- (a3);
 \draw (b1) -- (a4);
 \draw (b1) -- (a5);
 \draw (b1) -- (a6);
 \draw (b1) -- (a7);
 \draw (b1) -- (a8);
 \draw (b1) -- (a9);
 \draw (b1) -- (a10);
  
 \draw (b2) -- (a2);
 \draw (b2) -- (a2);
 \draw (b2) -- (a3);
 \draw (b2) -- (a4);
 \draw (b2) -- (a5);
 \draw (b2) -- (a6);
 \draw (b2) -- (a7);
 \draw (b2) -- (a8);
 \draw (b2) -- (a9);
 \draw (b2) -- (a10);

 \draw[very thick] (b2) -- (b1);
 
 \draw (b3) -- (b4);
 \draw (b1) -- (b4);
 \draw (b3) -- (b2);

    \node [dark] (c1) at (9.5,-3){} ;  
  \node [dark] (c2) at (13,-4)  {}; 
  \node [dark] (c3) at (16.5,-3)  {};  
  \node [dark] (c4) at (11.5,-1) {};  
  \node [dark] (c5) at (14.5,-1)  {};  
  \node [dark] (c6) at (11.5,1)  {};  
  \node [dark] (c7) at (14.5,1)  {};  
  \node [dark] (c8) at (9.5,3){};
  \node [dark] (c9) at (13,4) {};
  \node [dark] (c10) at (16.5,3) {};
  
  \node [hollow] (d1) at (7.5,0) {};
   \node [hollow] (d2) at (18.5,0) {};

   \node [hollow] (d3) at (7.5,-5) {};
   \node [hollow] (d4) at (18.5,-5) {};

    \draw (d1) -- (c2);
 \draw (d1) -- (c2);
 \draw (d1) -- (c3);
 \draw (d1) -- (c4);
 \draw (d1) -- (c5);
 \draw (d1) -- (c6);
 \draw (d1) -- (c7);
 \draw (d1) -- (c8);
 \draw (d1) -- (c9);
 \draw (d1) -- (c10);
  
 \draw (d2) -- (c2);
 \draw (d2) -- (c2);
 \draw (d2) -- (c3);
 \draw (d2) -- (c4);
 \draw (d2) -- (c5);
 \draw (d2) -- (c6);
 \draw (d2) -- (c7);
 \draw (d2) -- (c8);
 \draw (d2) -- (c9);
 \draw (d2) -- (c10);
  
   \draw[very thick] (d2) -- (d1);
    \draw (d3) -- (d4);
     \draw (d1) -- (d4);
      \draw (d3) -- (d2);
  
  \draw[red, very thick] (c1) -- (c2);
  \draw[red, very thick] (c2) -- (c3);  
  \draw[red, very thick] (c2) -- (c4);  
  \draw[red, very thick] (c4) -- (c6);  
  \draw[red, very thick] (c3) -- (c7);  
  \draw[red, very thick] (c6) -- (c8);  
  \draw[red, very thick] (c5) -- (c7);
  \draw[red, very thick] (c2) -- (c5);  
  \draw[red, very thick] (c8) -- (c4);  
  \draw[red, very thick] (c9) -- (c6);  
  \draw[red, very thick] (c3) -- (c4);  
  \draw[red, very thick] (c9) -- (c7);  
  \draw[red, very thick] (c5) -- (c1); 
  \draw[red, very thick] (c10) -- (c3);  
  \draw[red, very thick] (c1) -- (c8);  
  \draw[red, very thick] (c1) -- (c6);  
  \draw[red, very thick] (c10) -- (c5);  
  \draw[red, very thick] (c10) -- (c7);  
  \draw[red, very thick] (c8) -- (c9);
  \draw[red, very thick] (c10) -- (c9);

\end{tikzpicture}  
\end{center}
\caption{Non-isomorphic non-regular $A_\alpha$-cospectral graphs} \label{alcoregnshj}
\end{figure}
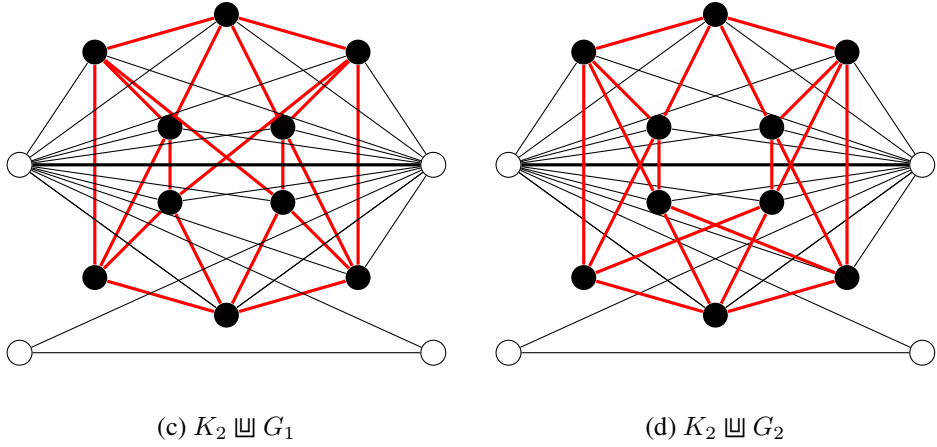

Using Corollary \ref{cosnshj}, an example of non-isomorphic non-regular $A_\alpha$-cospectral graph is shown in Figure \ref{alcoregnshj}, where $G_1$ and $G_2$ are graphs shown in Figure \ref{acoreg}. 

\begin{corollary}

    \begin{enumerate}
    \item Let $G_1$ and $G_2$ be two $A$-cospectral regular graphs and $H$ be an arbitrary graph. Then the graphs $G_1\doublesqcap H$ and $G_2\doublesqcap H$ are $A_\alpha$-cospectral.
    \item Let $H_1$ and $H_2$ be two $A_\alpha$-cospectral graphs with $\Gamma_{A_\alpha(H_1)}(\lambda-\alpha n_1) = \Gamma_{A_\alpha(H_2)}(\lambda-\alpha n_1)$ for $\alpha \in [0, 1]$. If $G$ is a regular graph, then the graphs $G\doublesqcap H_1$ and $G\doublesqcap H_2$ are $A_\alpha$-cospectral.
\end{enumerate}
\end{corollary}

\subsection{Central joins}
In this section, we derive the $A_\alpha$-characteristic polynomial of $G_1 \veedot G_2$ and $G_1 \veebar G_2$. 

\begin{proposition}\label{prop1}
Let $G_1$ be an $r_1$ regular graph with $n_1$ vertices and $\displaystyle m_1$ edges and $G_2$ be an arbitrary graph on $n_2$ vertices. Then the $A_\alpha$-characteristic polynomial of $G_1\veedot G_2$ is
\begin{align*}
\phi(A_\alpha(G_1\veedot G_2),\lambda)=&(\lambda-2\alpha)^{m_1-n_1}\prod_{i=1}^{n_2}\Big(\lambda-n_1\alpha-\lambda_i(A_{\alpha_2}))\Big)\\
&\prod_{j=2}^{n_1}\Big((\lambda-2\alpha)(\lambda-\alpha(n_1+n_2+\lambda_j(A_1))+1+\lambda_j(A_1))-(1-\alpha)^2(r_1+\lambda_j(A_1))\Big)\\
&\hspace{-1cm}\Big((\lambda-2\alpha)(\lambda-n_1(1+(1-\alpha)\Gamma_{A_{\alpha_2}}(\lambda-\alpha n_1))-\alpha n_2+(1-\alpha)r_1+1)-2(1-\alpha)^2r_1\Big).
\end{align*}

\end{proposition}
\begin{proof}
The $A_\alpha$ matrix of central vertex join of two graphs $G_1$ and $G_2$ is of the form
\begin{align*}
A_\alpha(G_1\veedot G_2)&=\begin{bmatrix}
\alpha(n_1+n_2-1)I_{n_1}+(1-\alpha)\overline{A_1} & (1-\alpha)R_1 & (1-\alpha)J_{n_1\times n_2}\\
(1-\alpha)R_1^T & 2\alpha I_{m_1} & 0_{m_1\times n_2}\\
(1-\alpha)J_{n_2\times n_1} & 0_{n_2\times m_1} & \alpha n_1I_{n_2}+A_{\alpha_2}
\end{bmatrix}\\
\intertext{Then}
\phi(A_\alpha(G_1\veedot& G_2),\lambda)=\begin{vmatrix}
(\lambda-\alpha(n_1+n_2-1))I_{n_1}-(1-\alpha)\overline{A_1} & -(1-\alpha)R_1 & -(1-\alpha)J_{n_1\times n_2}\\
-(1-\alpha)R_1^T & (\lambda-2\alpha) I_{m_1} & 0_{m_1\times n_2}\\
-(1-\alpha)J_{n_2\times n_1} & 0_{n_2\times m_1} & (\lambda-\alpha n_1)I_{n_2}-A_{\alpha_2}
\end{vmatrix}\\
&=|(\lambda-\alpha n_1)I_{n_2}-A_{\alpha_2}|\det{S}\\
&=\prod_{i=1}^{n_2}\Big(\lambda-\alpha n_1-\lambda_i(A_{\alpha_2})\Big)\det{S},\\
 \text { where } S&=\left[\begin{array}{ccc}(\lambda-\alpha(n_1+n_2-1)) I_{n_{1}}-(1-\alpha)(J_{n_{1}}-I_{n_{1}}-A_1) & -(1-\alpha)R_1 \\ -(1-\alpha)R_1^{T} & (\lambda-2\alpha) I_{m 1}\end{array}\right]\\ 
 &\qquad -\left[\begin{array}{c}-(1-\alpha)J_{n_{1} \times n_{2}} \\ O_{m_{1} \times n_{2}}\end{array}\right]\left((\lambda-\alpha n_1) I_{n 2}-A_{\alpha_2}\right)^{-1}\left[-(1-\alpha)J_{n_{2} \times n_{1}} \quad O_{n_{2} \times m_{1}}\right] \\
 &=\left[\begin{array}{cc}(\lambda-\alpha(n_1+n_2-1)) I_{n_{1}}-(1-\alpha)(J_{n_{1}}-I_{n_{1}}-A_1)-(1-\alpha)^2\Gamma_{A_{\alpha_2}}(\lambda-\alpha n_1)J_{n_1}  \\ -(1-\alpha)R_1^{T} \end{array}\right.\\
 &\hspace{9cm}\left. \begin{array}{c}
    -(1-\alpha)R_1    \\
       (\lambda-2\alpha) I_{m 1}
 \end{array}\right]\\
 \end{align*}
 \begin{align*}
 \det(S)&=(\lambda-2\alpha)^{m_1}\bigg|(\lambda-\alpha(n_1+n_2-1))I-(1-\alpha)(J-I-A_1)\\
 &\left.\hspace{5cm}-(1-\alpha)^2\Gamma_{A_{\alpha_2}}(\lambda-\alpha n_1)J_{n_1}-\frac{(1-\alpha)^2R_1R_1^T}{\lambda-2\alpha}\right|\\
 &=(\lambda-2\alpha)^{m_1}\bigg|(\lambda-\alpha(n_1+n_2-1))I-(1-\alpha)(P(A_1)-I-A_1)\\
 &\left.\hspace{5cm}-(1-\alpha)^2\Gamma_{A_{\alpha_2}}(\lambda-\alpha n_1)P(A_1)-\frac{(1-\alpha)^2(A_1 + r_1I)}{\lambda-2\alpha}\right|\\
 &=(\lambda-2\alpha)^{m_1}\prod_{i=1}^{n_1}\Big(\lambda-\alpha(n_1+n_2-1)-(1-\alpha)(P(\lambda_i(A_1))-1-\lambda_i(A_1))\\
 &\hspace{4cm}-(1-\alpha)^2\Gamma_{A_{\alpha_2}}(\lambda-\alpha n_1)P(\lambda_i(A_1))-\frac{(1-\alpha)^2(\lambda_i(A_1)+r_1)}{\lambda-2\alpha}\Big)\\
  &\hspace{-1cm}=(\lambda-2\alpha)^{m_1-n_1}\Big((\lambda-2\alpha)(\lambda-\alpha(n_2+r_1)-n_1+1+r_1-n_1(1-\alpha)^2\Gamma_{A_{\alpha_2}}(\lambda-\alpha n_1))-2r_1(1-\alpha)^2\Big)\\
 &\prod_{i=2}^{n_1}\Big((\lambda-2\alpha)(\lambda-\alpha(n_1+n_2)+(1-\alpha)\lambda_i(A_1)+1)-(1-\alpha)^2(\lambda_i(A_1)+r_1)\Big).\\
\end{align*}
Thus,
\begin{align*}
\phi(A_\alpha(G_1\veedot G_2),\lambda)=&(\lambda-2\alpha)^{m_1-n_1}\prod_{i=1}^{n_2}\Big(\lambda-n_1\alpha-\lambda_i(A_{\alpha_2})\Big)\\
&\prod_{j=2}^{n_1}\Big((\lambda-2\alpha)(\lambda-\alpha(n_1+n_2+\lambda_j(A_1))+1+\lambda_j(A_1))-(1-\alpha)^2(r_1+\lambda_j(A_1))\Big)\\
&\hspace{-1cm}\Big((\lambda-2\alpha)(\lambda-n_1(1+(1-\alpha)^2\Gamma_{A_{\alpha_2}}(\lambda-\alpha n_1))-\alpha n_2+(1-\alpha)r_1+1)-2(1-\alpha)^2r_1\Big).
\end{align*}

\end{proof}
As a consequence of Proposition \ref{prop1}, we obtain the following corollary.
\begin{corollary}
The $A_\alpha$-characteristic polynomial of $G_1\veedot G_2$ for two regular graphs $G_1$ and $G_2$ is
\begin{align*}
\phi(A_\alpha(G_1\veedot G_2),\lambda)=&(\lambda-2\alpha)^{m_1-n_1}\prod_{i=2}^{n_2}\Big(\lambda-n_1\alpha-\lambda_i(A_{\alpha_2})\Big)\\
&\prod_{j=2}^{n_1}\Big((\lambda-2\alpha)(\lambda-\alpha(n_1+n_2+\lambda_j(A_1))+1+\lambda_j(A_1))-(1-\alpha)^2(r_1+\lambda_j(A_1))\Big)\\
&\Big((\lambda-2\alpha)((\lambda-\alpha n_1-r_2)(\lambda-\alpha n_2+(1-\alpha)r_1+1)-n_1(\lambda-\alpha n_1-r_2+(1-\alpha)n_2))\\
&\hspace{8.5cm}-2r_1(1-\alpha)^2(\lambda-\alpha n_1-r_2)\Big).
\end{align*}
\end{corollary}

Now, in the following corollary, we obtain the $A_\alpha$-eigenvalues of $G_1\veedot G_2$, where $G_1$ and $G_2$ are regular graphs.

\begin{corollary}
Let $G_i$ be an $r_i$-regular graph with $n_i$ vertices and $m_i$ edges for $i = 1, 2$. Then the $A_\alpha$-spectrum of $G_1\veedot G_2$ consists of:
\begin{enumerate}
    \item $2\alpha$ repeated $m_1-n_1$ times,
    \item $n_1\alpha+\lambda_i(A_{\alpha_2})$, $i=2,3,\dots,n_2$,
    \item $\alpha +\displaystyle\frac{\alpha(n_1+n_2)-\lambda_i(A_1)(1-\alpha)-1}{2}$ \\ 
    $\pm\frac{\sqrt{(\lambda_i(A_1)(1-\alpha) + 1)^2+4(\lambda_i(A_1)(1-\alpha)+r_1)+\alpha^2((n_1 + n_2-2)^2+2(2 r_1 + 1 + \lambda_i(A_1)(n_1 + n_2)))+2\alpha(2 -(1+\lambda_i(A_1))(n_1 + n_2) - 4r_1)}}{2}$, $i=2,3,\dots,n_1$ and 
    \item three roots of the equation $(\lambda-2\alpha)((\lambda-\alpha n_1-r_2)(\lambda-\alpha n_2+(1-\alpha)r_1+1)-n_1(\lambda-\alpha n_1-r_2+(1-\alpha)n_2))-2r_1(1-\alpha)^2(\lambda-\alpha n_1-r_2)=0$.
\end{enumerate}
\end{corollary}

In the following corollary, we obtain the $A_\alpha$-eigenvalues of $G_1\veedot G_2$, where $G_1$ is a regular graph and $G_2$ is a non-regular graph.

\begin{corollary}
 Let $G_1$ be an $r_1$-regular graph on $n_1$ vertices and $m_1$ edges. Then the $A_\alpha$-spectrum of $G_1\veedot K_{p,q}$ is
 \begin{enumerate}
     \item $2\alpha$ repeated $m_1-n_1$ times,
     \item $\alpha(n_1+p)$ repeated $q-1$ times,
     \item $\alpha(n_1+q)$ repeated $p-1$ times,
      \item $\alpha +\displaystyle\frac{\alpha(n_1+p+q)-\lambda_i(A_1)(1-\alpha)-1}{2}$ \\ 
    $\pm\frac{\sqrt{(\lambda_i(A_1)(1-\alpha) + 1)^2+4(\lambda_i(A_1)(1-\alpha)+r_1)+\alpha^2((n_1 + p+q-2)^2+2(2 r_1 + 1 + \lambda_i(A_1)(n_1 + p+q)))+2\alpha(2 -(1+\lambda_i(A_1))(n_1 + p+q) - 4r_1)}}{2}$, $i=2,3,\dots,n_1$ and 
    \item three roots of the equation $(\lambda-2\alpha)\Big(((\lambda-\alpha n_1)^2-\alpha(p+q)(\lambda-\alpha n_1)+(2\alpha-1)pq)(\lambda-n_1-\alpha(p+q)+(1-\alpha)r_1+1)-n_1(1-\alpha)((p+q)(\lambda-\alpha n_1)-\alpha(p+q)^2+2pq)\Big)-2r_1(1-\alpha)^2((\lambda-\alpha n_1)^2-\alpha(p+q)(\lambda-\alpha n_1)+(2\alpha-1)pq)=0$.
 \end{enumerate}
\end{corollary}

The following corollary presents new pairs of non-isomorphic non-regular $A_\alpha$-cospectral graphs.

\begin{corollary}

    \begin{enumerate}
    \item Let $G_1$ and $G_2$ be two $A$-cospectral regular graphs and $H$ be an arbitrary graph. Then the graphs $G_1\veedot H$ and $G_2\veedot H$ are $A_\alpha$-cospectral.
    \item Let $H_1$ and $H_2$ be two $A_\alpha$-cospectral graphs with $\Gamma_{A_\alpha(H_1)}(x) = \Gamma_{A_\alpha(H_2)}(x)$ for $\alpha \in [0, 1]$. If $G$ is a regular graph, then the graphs $G\veedot H_1$ and $G\veedot H_2$ are $A_\alpha$-cospectral.
\end{enumerate}
\end{corollary}

\begin{proposition}\label{cej}
Let $G_i$ be an $r_i$ regular graph with $n_i$ vertices and $\displaystyle m_i$ edges. Then the $A_\alpha$-characteristic polynomial of $G_1\veebar G_2$ is $\phi(A_\alpha(G_1\veebar G_2),\lambda)=$
\begin{align*}
&(\lambda-\alpha(2+n_2))^{m_1-n_1-1}\\
&\big(((\lambda-\alpha(2+n_2))(\lambda-\alpha m_1-r_2)-m_1n_2(1-\alpha)^2)((\lambda-\alpha(2+n_2))(\lambda-n_1+1+(1-\alpha)r_1)-2r_1(1-\alpha)^2)\\
 &\hspace{13.5cm}-n_1n_2r_1^2(1-\alpha)^4\big)\\
 &\prod_{i=2}^{n_2}(\lambda-\alpha(m_1+r_2)-(1-\alpha)\lambda_i(A_2))\\
 &\prod_{i=2}^{n_1}\left((\lambda-\alpha(2+n_2))(\lambda-\alpha(n_1-1)+(1-\alpha)(1+\lambda_i(A_1)))-(1-\alpha)^2(\lambda_i(A_1)+r_1)\right)
 \end{align*}
\end{proposition}

\begin{proof}
The $A_\alpha$ matrix of central edge join of two graphs $G_1$ and $G_2$ is of the form
\begin{align*}
A_\alpha(G_1\veebar G_2)&=\begin{bmatrix}
\alpha(n_1-1)I_{n_1}+(1-\alpha)\overline{A_1} & (1-\alpha)R_1 & O_{n_1\times n_2}\\
(1-\alpha)R_1^T & \alpha(2+n_2) I_{m_1} & (1-\alpha)J_{m_1\times n_2}\\
O_{n_2\times n_1} & (1-\alpha)J_{n_2\times m_1} & \alpha m_1I_{n_2}+A_{\alpha_2}
\end{bmatrix}\\
\intertext{Then}
\phi(A_\alpha(G_1\veebar& G_2),\lambda)=\begin{vmatrix}
(\lambda-\alpha(n_1-1))I_{n_1}-(1-\alpha)\overline{A_1} & -(1-\alpha)R_1 & O_{n_1\times n_2}\\
-(1-\alpha)R_1^T & (\lambda-\alpha(2+n_2)) I_{m_1} & -(1-\alpha)J_{m_1\times n_2}\\
O_{n_2\times n_1} & -(1-\alpha)J_{n_2\times m_1} & (\lambda-\alpha m_1)I_{n_2}-A_{\alpha_2}
\end{vmatrix}\\
&=|(\lambda-\alpha m_1)I_{n_2}-A_{\alpha_2}|\det{S}\\
&=\prod_{i=1}^{n_2}\Big(\lambda-\alpha m_1-\lambda_i(A_{\alpha_2})\Big)\det{S},\\
 \text { where } S&=\left[\begin{array}{ccc}(\lambda-\alpha(n_1-1))I_{n_1}-(1-\alpha)\overline{A_1} & -(1-\alpha)R_1 \\ -(1-\alpha)R_1^{T} & (\lambda-\alpha(2+n_2)) I_{m_1}-(1-\alpha)^2\Gamma_{A_{\alpha_2}}(\lambda-\alpha m_1)J_{m_1}
 \end{array}\right]\\ 
 \end{align*}
 \begin{align*}
 \det(S)=&\left|(\lambda-\alpha(2+n_2)) I_{m_1}-(1-\alpha)^2\Gamma_{A_{\alpha_2}}(\lambda-\alpha m_1)J_{m_1}\right|\\
 &\bigg|(\lambda-\alpha(n_1-1))I_{n_1}-(1-\alpha)\overline{A_1}-(1-\alpha)^2R_1((\lambda-\alpha(2+n_2)) I_{m_1}-(1-\alpha)^2\Gamma_{A_{\alpha_2}}(\lambda-\alpha m_1)J_{m_1})^{-1}R_1^T\bigg|\\
 \intertext{using Lemma \ref{inverse}}
 \det(S)=&\left|(\lambda-\alpha(2+n_2)) I_{m_1}-(1-\alpha)^2\Gamma_{A_{\alpha_2}}(\lambda-\alpha m_1)J_{m_1}\right|\\
 &\bigg|(\lambda-\alpha(n_1-1))I_{n_1}-(1-\alpha)\overline{A_1}-\frac{(1-\alpha)^2R_1R_1^T}{\lambda-\alpha(2+n_2)}\\
 &\hspace{6cm}-\frac{(1-\alpha)^4\Gamma_{A_{\alpha_2}}(\lambda-\alpha m_1)r_1^2J_{n_1}}{(\lambda-\alpha(2+n_2))(\lambda-\alpha(2+n_2)-m_1(1-\alpha)^2\Gamma_{A_{\alpha_2}}(\lambda-\alpha m_1))}\bigg|\\
 =&\left|(\lambda-\alpha(2+n_2)) I_{m_1}-(1-\alpha)^2\Gamma_{A_{\alpha_2}}(\lambda-\alpha m_1)J_{m_1}\right|\\
 &\bigg|(\lambda-\alpha(n_1-1))I_{n_1}-(1-\alpha)\overline{A_1}-\frac{(1-\alpha)^2(A_1+r_1I_{n_1})}{\lambda-\alpha(2+n_2)}\\
 &\hspace{6cm}-\frac{(1-\alpha)^4\Gamma_{A_{\alpha_2}}(\lambda-\alpha m_1)r_1^2P(A_1)}{(\lambda-\alpha(2+n_2))(\lambda-\alpha(2+n_2)-m_1(1-\alpha)^2\Gamma_{A_{\alpha_2}}(\lambda-\alpha m_1))}\bigg|\\
 =&\frac{(\lambda-\alpha(2+n_2))^{m_1-1}\left((\lambda-\alpha m_1-r_2)(\lambda-\alpha(2+n_2))-(1-\alpha)^2m_1n_2\right)}{\lambda-\alpha m_1-r_2}\\
 &\frac{1}{(\lambda-\alpha(2+n_2))^{n_1}((\lambda-\alpha(2+n_2))(\lambda-\alpha m_1-r_2)-m_1n_2(1-\alpha)^2)}\\
 &\big(((\lambda-\alpha(2+n_2))(\lambda-\alpha m_1-r_2)-m_1n_2(1-\alpha)^2)((\lambda-\alpha(2+n_2))(\lambda-n_1+1+(1-\alpha)r_1)-2r_1(1-\alpha)^2)\\
 &\hspace{13.5cm}-n_1n_2r_1^2(1-\alpha)^4\big)\\
 &\prod_{i=2}^{n_1}\left((\lambda-\alpha(2+n_2))(\lambda-\alpha(n_1-1)+(1-\alpha)(1+\lambda_i(A_1)))-(1-\alpha)^2(\lambda_i(A_1)+r_1)\right) 
\end{align*}

Thus, $\phi(A_\alpha(G_1\veebar G_2),\lambda)=$
\begin{align*}
&(\lambda-\alpha(2+n_2))^{m_1-n_1-1}\\
&\big(((\lambda-\alpha(2+n_2))(\lambda-\alpha m_1-r_2)-m_1n_2(1-\alpha)^2)((\lambda-\alpha(2+n_2))(\lambda-n_1+1+(1-\alpha)r_1)-2r_1(1-\alpha)^2)\\
 &\hspace{13.5cm}-n_1n_2r_1^2(1-\alpha)^4\big)\\
 &\prod_{i=2}^{n_2}(\lambda-\alpha(m_1+r_2)-(1-\alpha)\lambda_i(A_2))\\
 &\prod_{i=2}^{n_1}\left((\lambda-\alpha(2+n_2))(\lambda-\alpha(n_1-1)+(1-\alpha)(1+\lambda_i(A_1)))-(1-\alpha)^2(\lambda_i(A_1)+r_1)\right)
\end{align*}

Now, in the following corollary, we obtain the $A_\alpha $-eigenvalues of $G_1 \veebar G_2$ where $G_i$'s are $r_i$-regular graphs.
\begin{corollary}
Let $G_i$ be an $r_i$-regular graph with $n_i$ vertices and $m_i$ edges for $i = 1, 2$. Then the $A_\alpha$-spectrum of $G_1\veebar G_2$ consists of:
\begin{enumerate}
    \item $\alpha(2+n_2)$ repeated $m_1-n_1-1$ times,
    \item $m_1\alpha+\lambda_i(A_{\alpha_2})$, $i=2,3,\dots,n_2$,
    \item two roots of the equation $(\lambda-\alpha(2+n_2))(\lambda-\alpha(n_1-1)+(1-\alpha)(1+\lambda_i(A_1)))-(1-\alpha)^2(\lambda_i(A_1)+r_1)=0$ for $i=2,3,\dots,n_2$ and 
    \item four roots of the equation $((\lambda-\alpha(2+n_2))(\lambda-\alpha m_1-r_2)-m_1n_2(1-\alpha)^2)((\lambda-\alpha(2+n_2))(\lambda-n_1+1+(1-\alpha)r_1)-2r_1(1-\alpha)^2)-n_1n_2r_1^2(1-\alpha)^4=0$.
\end{enumerate}
\end{corollary}

The subsequent corollary introduces new pairs of $A_\alpha$-cospectral graphs that are not regular.

\begin{corollary}

    \begin{enumerate}
    \item Let $G_1$ and $G_2$ be two $A$-cospectral regular graphs and $H$ be any regular graph. Then the graphs $G_1\veebar H$ and $G_2\veebar H$ are $A_\alpha$-cospectral.
    \item Let $H_1$ and $H_2$ be two $A_\alpha$-cospectral regular graphs. If $G$ is a regular graph, then the graphs $G\veebar H_1$ and $G\veebar H_2$ are $A_\alpha$-cospectral.
    \item Let $G_1$ and $G_2$ be $A_\alpha$-cospectral regular graphs, and $H_1$ and $H_2$ be another $A_\alpha$-cospectral regular graphs. Then $G_1\veebar H_1$ and $G_2\veebar H_2$ are $A_\alpha$-cospectral.
\end{enumerate}
\end{corollary}

\end{proof}
\subsection{Duplicate join}
In this section, we derive the $A_\alpha$-characteristic polynomial of $G_1 \bowtie G_2$ when $G_1$ and $G_2$ are arbitrary graphs.
\begin{proposition}\label{double}
    Let $G_i$ be a graph on $n_i$ vertices for $i=1,2$. Then

\begin{align*}
\phi(A_\alpha  \left(G_1 \bowtie G_2\right),\lambda) =&\left|(\lambda-\alpha n_1)I-A_{\alpha_2}\right| \left|\lambda I-\alpha D_1\right| \\
& \left|(\lambda-\alpha n_2)I-\alpha D_1-(1-\alpha)^2 A_1(\lambda I-\alpha D_1)^{-1}A_1\right| \\
& \left[1-\Gamma_{A_{\alpha_2}}\left(\lambda-\alpha n_1\right) \Gamma_{\alpha D_1+(1-\alpha)^2A_1(\lambda I-\alpha D_1)^{-1}A_1}\left(\lambda-\alpha n_2\right)\right] .
\end{align*}

\end{proposition}

\begin{proof}
   With a suitable ordering of the vertices of $G_1 \bowtie G_2$, we get

\begin{align*}
A_\alpha\left(G_1 \bowtie G_2\right)=&\left(\begin{array}{ccc}
\alpha n_2 I_{n_1}+\alpha D_1 & (1-\alpha)A_1 & (1-\alpha)J_{n_1 \times n_2} \\
(1-\alpha)A_1 & \alpha D_1 & O_{n_1 \times n_2} \\
(1-\alpha)J_{n_2 \times n_1} & O_{n_2 \times n_1} & \alpha n_1 I_{n_2}+A_{\alpha_2}
\end{array}\right).
\end{align*}
Then the results follow from Lemma \ref{schurmy}.

\end{proof}

Now, in the following propositions, we obtain the $A_\alpha $-eigenvalues of $G_1 \bowtie G_2$ where $G_i$'s are $r_i$-regular graphs.

\begin{proposition}
    Let $G_i$ be $r_i$-regular graph with $n_i$ vertices for $i=1,2$. Then the $A_\alpha$-spectrum of $G_1 \bowtie G_2$ consists of:
    \begin{enumerate}
        \item $\alpha(n_1+r_2)+(1-\alpha)\lambda_i(A_2)$ for each $i=2,3,\dots,n_2$,
        \item two roots of the equation $$(\lambda-\alpha r_1)(\lambda-\alpha(n_2+r_1))-(1-\alpha)^2\lambda_i(A_1)=0, \text{ for each $i=2,3\dots,n_1$,}$$
        \item three roots of the equation
        $$(\lambda-\alpha n_1-r_2)((\lambda-\alpha r_1)(\lambda-\alpha(n_2+r_1))-(1-\alpha)^2r_1^2)-n_1n_2(\lambda-\alpha r_1)=0.$$
    \end{enumerate}
\end{proposition}

The following corollary presents new pairs of non-isomorphic non-regular $A_\alpha$-cospectral graphs.

\begin{corollary}

    \begin{enumerate}
    \item Let $G_1$ and $G_2$ be two $A_\alpha$-cospectral regular graphs and $H$ be an arbitrary graph. Then the graphs $G_1\bowtie H$ and $G_2\bowtie H$ are $A_\alpha$-cospectral.
    \item Let $H_1$ and $H_2$ be two $A_\alpha$-cospectral graphs with $\Gamma_{A_\alpha(H_1)}(\lambda-\alpha n_1) = \Gamma_{A_\alpha(H_2)}(\lambda-\alpha n_1)$ for $\alpha \in [0, 1]$. If $G$ is a regular graph, then the graphs $G\bowtie H_1$ and $G\bowtie H_2$ are $A_\alpha$-cospectral.
\end{enumerate}
\end{corollary}

\begin{figure}[H]
\begin{center}
\begin{tikzpicture}[scale=.5,auto=center] 
 \tikzset{dark/.style={circle,fill=black}}
 \tikzset{hollow/.style={circle,draw=black}}
 \tikzset{white/.style={circle,draw=white}}
    \node [dark] (a1) at (-3.5,-3){} ;  
  \node [dark] (a2) at (0,-4)  {}; 
  \node [dark] (a3) at (3.5,-3)  {};  
  \node [dark] (a4) at (-1.5,-1) {};  
  \node [dark] (a5) at (1.5,-1)  {};  
  \node [dark] (a6) at (-1.5,1)  {};  
  \node [dark] (a7) at (1.5,1)  {};  
  \node [dark] (a8) at (-3.5,3){};
  \node [dark] (a9) at (0,4) {};
  \node [dark] (a10) at (3.5,3) {};

  \node [hollow] (b1) at (-5.5,0){} ;  
   \node [hollow] (b2) at (5.5,0){} ;  

   \node [hollow] (b3) at (-5.5,-5){} ;  
   \node [hollow] (b4) at (5.5,-5){} ;  
  
  \node  at (0,-7) {(c) $K_2\bowtie G_1$};
   \node  at (13,-7) {(d) $K_2\bowtie G_2$};
  
  \draw[red, very thick] (a1) -- (a2);
  \draw[red, very thick] (a2) -- (a3);  
  \draw[red, very thick] (a2) -- (a4);  
  \draw[red, very thick] (a4) -- (a6);  
  \draw[red, very thick] (a3) -- (a7);  
  \draw[red, very thick] (a6) -- (a8);  
  \draw[red, very thick] (a5) -- (a7);
  \draw[red, very thick] (a2) -- (a5);  
  \draw[red, very thick] (a8) -- (a5);  
  \draw[red, very thick] (a9) -- (a6);  
  \draw[red, very thick] (a3) -- (a5);  
  \draw[red, very thick] (a9) -- (a7);  
  \draw[red, very thick] (a4) -- (a1); 
  \draw[red, very thick] (a10) -- (a3);  
  \draw[red, very thick] (a1) -- (a8);  
  \draw[red, very thick] (a1) -- (a6);  
  \draw[red, very thick] (a10) -- (a4);  
  \draw[red, very thick] (a10) -- (a7);  
  \draw[red, very thick] (a8) -- (a9);
  \draw[red, very thick] (a10) -- (a9);
  
  \draw (b1) -- (a2);
 \draw (b1) -- (a2);
 \draw (b1) -- (a3);
 \draw (b1) -- (a4);
 \draw (b1) -- (a5);
 \draw (b1) -- (a6);
 \draw (b1) -- (a7);
 \draw (b1) -- (a8);
 \draw (b1) -- (a9);
 \draw (b1) -- (a10);
  
 \draw (b2) -- (a2);
 \draw (b2) -- (a2);
 \draw (b2) -- (a3);
 \draw (b2) -- (a4);
 \draw (b2) -- (a5);
 \draw (b2) -- (a6);
 \draw (b2) -- (a7);
 \draw (b2) -- (a8);
 \draw (b2) -- (a9);
 \draw (b2) -- (a10);

 
 \draw (b1) -- (b4);
 \draw (b3) -- (b2);

    \node [dark] (c1) at (9.5,-3){} ;  
  \node [dark] (c2) at (13,-4)  {}; 
  \node [dark] (c3) at (16.5,-3)  {};  
  \node [dark] (c4) at (11.5,-1) {};  
  \node [dark] (c5) at (14.5,-1)  {};  
  \node [dark] (c6) at (11.5,1)  {};  
  \node [dark] (c7) at (14.5,1)  {};  
  \node [dark] (c8) at (9.5,3){};
  \node [dark] (c9) at (13,4) {};
  \node [dark] (c10) at (16.5,3) {};
  
  \node [hollow] (d1) at (7.5,0) {};
   \node [hollow] (d2) at (18.5,0) {};

   \node [hollow] (d3) at (7.5,-5) {};
   \node [hollow] (d4) at (18.5,-5) {};

    \draw (d1) -- (c2);
 \draw (d1) -- (c2);
 \draw (d1) -- (c3);
 \draw (d1) -- (c4);
 \draw (d1) -- (c5);
 \draw (d1) -- (c6);
 \draw (d1) -- (c7);
 \draw (d1) -- (c8);
 \draw (d1) -- (c9);
 \draw (d1) -- (c10);
  
 \draw (d2) -- (c2);
 \draw (d2) -- (c2);
 \draw (d2) -- (c3);
 \draw (d2) -- (c4);
 \draw (d2) -- (c5);
 \draw (d2) -- (c6);
 \draw (d2) -- (c7);
 \draw (d2) -- (c8);
 \draw (d2) -- (c9);
 \draw (d2) -- (c10);
  
     \draw (d1) -- (d4);
      \draw (d3) -- (d2);
  
  \draw[red, very thick] (c1) -- (c2);
  \draw[red, very thick] (c2) -- (c3);  
  \draw[red, very thick] (c2) -- (c4);  
  \draw[red, very thick] (c4) -- (c6);  
  \draw[red, very thick] (c3) -- (c7);  
  \draw[red, very thick] (c6) -- (c8);  
  \draw[red, very thick] (c5) -- (c7);
  \draw[red, very thick] (c2) -- (c5);  
  \draw[red, very thick] (c8) -- (c4);  
  \draw[red, very thick] (c9) -- (c6);  
  \draw[red, very thick] (c3) -- (c4);  
  \draw[red, very thick] (c9) -- (c7);  
  \draw[red, very thick] (c5) -- (c1); 
  \draw[red, very thick] (c10) -- (c3);  
  \draw[red, very thick] (c1) -- (c8);  
  \draw[red, very thick] (c1) -- (c6);  
  \draw[red, very thick] (c10) -- (c5);  
  \draw[red, very thick] (c10) -- (c7);  
  \draw[red, very thick] (c8) -- (c9);
  \draw[red, very thick] (c10) -- (c9);

\end{tikzpicture}  
\end{center}
\caption{Non-isomorphic non-regular $A_\alpha$-cospectral graphs} \label{alcoregdj}
\end{figure}
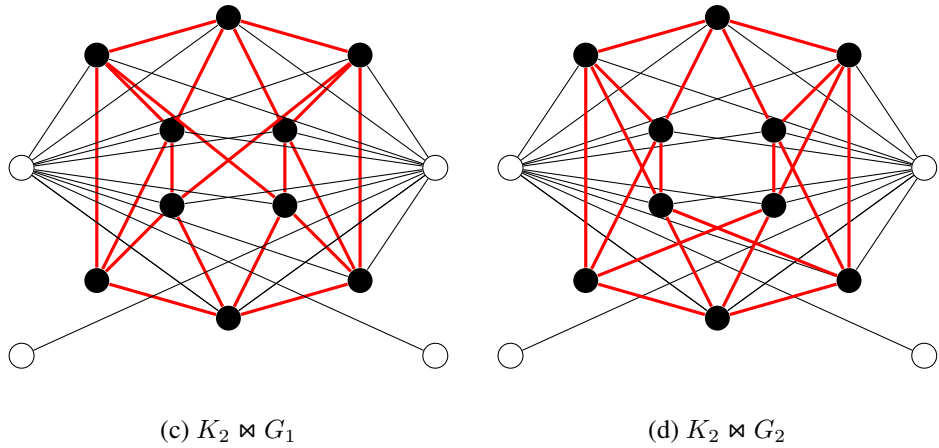

\section{Conclusion}
In this article, we have derived $A_\alpha$-characteristic polynomial and $A_\alpha$-spectrum for various join of two graphs such as neighbour and non-neighbour splitting join, neighbour and non-neighbour shadow join, central vertex and edge join, and duplicate join. Subsequently, we have generated infinitely many families of $A_\alpha$-cospectral graphs that are not regular.

\section*{Observations}
In this section, with the help of Matlab software, we study the $A_\alpha$-energy of some join of $K_2$ and $K_{p,p}$ as $\alpha$ varies. 

Finding the $A_\alpha$-energy of graphs is tedious when the graph has a complex structure. Here, we compute the $A_\alpha$-energy of some particular graphs and observe some relation with the $A_\alpha$-energy of the complete graph.

If $\bullet$ denotes any of the operation $\barwedge ,\doublebarwedge, \doublesqcup, \doublesqcap, \bowtie$, then all $K_2\bullet K_{p,p}$ has same number of vertices, that is $4+2p$ vertices.\\

\noindent\textbf{Observations:}
We observe the following from Table \ref{joinenergy}.
\begin{itemize}
    \item When $p\geq2$, $K_2\bullet K_{p,p}$ is $A_\alpha$-borderenergetic for $\alpha\geq0.4$.
    \item When $p=1$, $K_2\bullet K_{p,p}$ is $A_\alpha$-borderenergetic for $\alpha\geq0.5$.
\end{itemize}

\begin{table}[H]
    \centering
    \begin{tabular}{ccccccccccc}
        	&	0	&	0.1	&	0.2	&	0.3	&	0.4	&	0.5	&	0.6	&	0.7	&	0.8	&	0.9	\\
         \hline
         $K_6$	&	10	&	9	&	8	&	7	&	6	&	5	&	4	&	3	&	2	&	1	\\
$K_2\barwedge K_{1,1}$	&	8.0188	&	7.5096	&	7.2859	&	7.277	&	7.4041	&	7.6864	&	8.104	&	8.6173	&	9.1932	&	9.8105	\\
$K_2\doublebarwedge K_{1,1}$	&	6	&	6.2	&	6.4	&	6.6	&	6.8	&	7	&	7.6	&	8.2	&	8.8	&	9.4	\\
$K_2\doublesqcup K_{1,1}$	&	8.4721	&	7.9165	&	7.3936	&	6.9123	&	6.4824	&	6.1119	&	6.1495	&	6.3181	&	6.5432	&	6.8147	\\
$K_2\doublesqcap K_{1,1}$	&	8	&	7.4667	&	6.9333	&	6.4	&	5.8667	&	5.6667	&	5.6	&	6.2	&	6.8	&	7.4	\\
$K_2\bowtie K_{1,1}$	&	7.8064	&	7.1957	&	6.7487	&	6.3775	&	5.943	&	5.2649	&	5.7566	&	6.3996	&	6.9645	&	7.4926	\\
&&&&&&&&&&\\
$K_8$	&	14	&	12.6	&	11.2	&	9.8	&	8.4	&	7	&	5.6	&	4.2	&	2.8	&	1.4	\\
$K_2\barwedge K_{2,2}$	&	10.6885	&	10.8984	&	11.8064	&	12.9649	&	14.3742	&	15.9636	&	17.6646	&	19.436	&	21.2562	&	23.1132	\\
$K_2\doublebarwedge K_{2,2}$	&	8.7446	&	8.9131	&	9.0886	&	9.271	&	9.6604	&	10.4069	&	11.1604	&	11.921	&	12.6886	&	13.5131	\\
$K_2\doublesqcup K_{2,2}$	&	11.2078	&	10.4495	&	9.8095	&	9.2975	&	9.1462	&	9.4741	&	9.8805	&	10.3472	&	10.8616	&	11.4146	\\
$K_2\doublesqcap K_{2,2}$	&	10.7446	&	10.1131	&	9.4886	&	8.971	&	9.0604	&	9.6569	&	10.2604	&	10.871	&	11.4886	&	12.6131	\\
$K_2\bowtie K_{2,2}$	&	10.5941	&	9.9025	&	9.538	&	8.9296	&	9.2352	&	9.9558	&	10.5521	&	11.1091	&	11.6532	&	12.6963	\\
&&&&&&&&&&\\
$K_{10}$	&	18	&	16.2	&	14.4	&	12.6	&	10.8	&	9	&	7.2	&	5.4	&	3.6	&	1.8	\\
$K_2\barwedge K_{3,3}$	&	13.1213	&	14.9884	&	17.7574	&	20.9359	&	24.4066	&	28.0456	&	31.7874	&	35.6013	&	39.472	&	43.3902	\\
$K_2\doublebarwedge K_{3,3}$	&	11.2111	&	11.7506	&	12.3114	&	13.0542	&	14.3797	&	15.7282	&	17.0997	&	18.6542	&	20.9514	&	23.2706	\\
$K_2\doublesqcup K_{3,3}$	&	13.6696	&	12.7022	&	11.9752	&	11.6089	&	12.0728	&	12.6515	&	13.3155	&	14.0469	&	14.8346	&	15.6709	\\
$K_2\doublesqcap K_{3,3}$	&	13.2111	&	12.3906	&	11.5914	&	11.6142	&	12.4597	&	13.3282	&	14.2197	&	15.1342	&	16.0714	&	17.1506	\\
$K_2\bowtie K_{3,3}$	&	13.1056	&	12.3439	&	11.7362	&	11.8236	&	12.8467	&	13.7092	&	14.5408	&	15.3761	&	16.2292	&	17.2265	\\
&&&&&&&&&&\\
$K_{12}$	&	22	&	19.8	&	17.6	&	15.4	&	13.2	&	11	&	8.8	&	6.6	&	4.4	&	2.2	\\
$K_2\barwedge K_{4,4}$	&	15.435	&	19.7734	&	25.2475	&	31.2471	&	37.5426	&	44.0024	&	50.573	&	57.2305	&	63.9617	&	70.7577	\\
$K_2\doublebarwedge K_{4,4}$	&	13.544	&	14.8522	&	16.2	&	18.2895	&	20.6225	&	23	&	25.4225	&	28.9895	&	32.8	&	36.6522	\\
$K_2\doublesqcup K_{4,4}$	&	16	&	14.8257	&	14.0339	&	14.1923	&	14.8682	&	15.6729	&	16.578	&	17.5672	&	18.63	&	19.7581	\\
$K_2\doublesqcap K_{4,4}$	&	15.544	&	14.4856	&	13.7333	&	14.4895	&	15.5558	&	16.6667	&	17.8225	&	19.0228	&	20.2667	&	21.5522	\\
$K_2\bowtie K_{4,4}$	&	15.4785	&	14.5794	&	13.2076	&	14.9225	&	16.0284	&	17.0837	&	18.1542	&	19.2625	&	20.4175	&	21.6225	
    \end{tabular}
    \caption{$A_\alpha$-energy of $K_2\bullet K_{p,p}$}
    \label{joinenergy}
\end{table}
\bibliographystyle{unsrt}  
\bibliography{references}

\end{document}